\documentstyle[12pt]{article}
\begin{document}

\newcommand{\PSbox}[3]{\mbox{\rule{0in}{#3}\includegraphics{#1}\hspace{#2}}}
\newtheorem{theorem}{Theorem}[section]  
\newtheorem{definition}[theorem]{Definition}  
\newtheorem{example}[theorem]{Example}  
\newtheorem{lemma}[theorem]{Lemma}  
\newtheorem{proposition}[theorem]{Proposition}  
\newtheorem{corollary}[theorem]{Corollary}  
\newtheorem{remark}[theorem]{Remark}  
\newtheorem{conjecture}[theorem]{Conjecture}

\newcommand{\A}{{\cal A}} 
\newcommand{\B}{{\cal B}} 
\newcommand{\C}{{\cal C}} 
\newcommand{\E}{{\cal E}}
\newcommand{\F}{{\cal F}}
\newcommand{\Ho}{{\cal H}}
\newcommand{\M}{{\cal M}} 
\newcommand{\N}{{\cal N}} 
\newcommand{\V}{{\cal V}}
\newcommand{\Lo}{{\cal L}}
\newcommand{\X}{{\cal X}}
\newcommand{\al}{\alpha} 
\newcommand{\Ann}{\mbox{\rm Ann}}

\newcommand{\cn}{{\bf {\rm C}} 
\hspace{-.4em}      {\vrule height1.5ex width.08em depth-.04ex} 
\hspace{.3em}} 
\newcommand{\Hom}{\mbox{\rm Hom}}
\newcommand{\Der}{\mbox{\rm Der}(S)} 
\newcommand{\Char}{\mbox{\rm Char}}
\newcommand{\gl}{{g_{\lambda}}} 
\newcommand{\Gr}{\mbox{\rm Gr}} 
\newcommand{\ints}{{\sf Z}\hspace{-.36em}{\sf Z}} 
\newcommand{\la}{\lambda}
\newcommand{\om}{\omega}  
\newcommand{\Om}{\Omega}  
\newcommand{\Op}{\Omega^{p}}  
\newcommand{\OpA}{\Omega^{p}(\A)}

\newcommand{\p}{\partial}  
\newcommand{\PD}{\mbox{\rm pd\,}}
\newcommand{\Pl}{{\Phi_{\lambda}}}  
\newcommand{\pl}{{\phi_{\lambda}}}  
\newcommand{\Poin}{\mbox{\rm Poin}} 
\newcommand{\PoinO}{\Poin(\Omega^{*}(c\A); u, } 
\newcommand{\PoinGr}{\Poin(\Gr\Omega^{*}(\A); u, } 
\newcommand{\rk}{\mbox{\rm rank}}
\newcommand{\proof}{{\bf Proof.~}} 
\newcommand{\qed}{~~\mbox{$\Box$}} 

\newcommand{\ra}{{\rightarrow}}  
\newcommand{\rn}{{\rm I}\hspace{-.2em}{\rm R}}

\newcommand{\scn}{\scriptsize\cn} 
\newcommand{\stR}[1]{\stackrel{#1}{\longrightarrow}} 
\newcommand{\th}{\theta}

\newcommand{\we}{\wedge}
\newcommand{\ar}{\buildrel d_{\la}\over\rightarrow}

\title{Cohomology of the Orlik-Solomon algebras
and local systems}
\author  { 
{\sc Anatoly Libgober}\\
{\small\it Department of Mathematics, University of Illinois,
Chicago, Ill 60607}\\
{\sc Sergey Yuzvinsky }\\
{\small\it Department of Mathematics,   University of Oregon,
Eugene, OR 97403 } }

\date{\today} 
\maketitle 
\begin{abstract}
 The paper provides a combinatorial method to decide 
when the space of local systems with non vanishing 
first cohomology on the complement to 
an arrangement of lines in a complex projective plane 
has as an irreducible component a subgroup
of positive dimension.
Partial classification of arrangements having such a 
component of positive dimension and a comparison 
theorem for cohomology of Orlik-Solomon algebra 
and cohomology of local systems are given.
The methods are based on Vinberg-Kac classification of 
generalized Cartan matrices and study of pencils 
of algebraic curves 
defined by mentioned positive dimensional components.
 
\end{abstract}

\section{Introduction}
\bigskip

\par One of the central questions in the theory of hyperplane
arrangements is the
problem of expressing topological invariants of the complements to 
arrangements in terms of combinatorics.
Probably the first non trivial result in this direction is due to Arnold, 
Brieskorn and Orlik-Solomon who calculated the cohomology algebra of the
complement (referred below to as the Orlik-Solomon algebra) 
in terms of the intersection lattice of an
arrangement. The question to which extent the fundamental group 
of the complement can be described in combinatorial terms turned out 
to be rather difficult
(cf. \cite{Ko}, \cite{R}). In this paper we shall show how  
certain {\it invariants} of the fundamental group, namely its characteristic 
varieties containing the trivial character of the fundamental group,
can be calculated from the intersection lattice.
As for the most questions about fundamental 
groups, due to Lefschetz type theorems,  the study of these invariants 
sufficient to carry out for arrangements of lines. 
The arrangements for which
these invariants are non trivial have an interesting combinatorial
structure namely they admit certain partition of the collection of lines
into disjoint groups. We show that such decompositions can be detected either 
combinatorially, using methods going back to Vinberg and Kac  
and which were designed for the classification of Kac-Moody 
algebras, or algebro-geometrically, using pencils of planes curves 
whose singular members form the arrangement. 
Moreover we were able in some cases to give classification of 
arrangements for which these invariants are non trivial and it suggests 
that a complete classification may be a feasible task.
\par As we mentioned, the invariants of the fundamental group in question are 
the characteristic varieties. For an arrangement with the complement $M$
such a variety belongs to the torus of characters
of the fundamental group $Hom(\pi_1(M), {\bf C}^*)$. 
More precisely the characteristic variety of $M$ is formed by those 
characters of $\pi_1(M)$ for which the 
corresponding local system has non-vanishing cohomology in dimension one.   
It follows from a theorem of Arapura \cite{A} 
that these subvarieties are unions
of translated subtori and we are interested in subtori which 
contain the trivial character (i.e. for which no translation is needed).
\par The starting point in our study 
is the connection between the following two kinds of
cohomology that has been known for some time. The
first kind is the cohomology $H^*(A,a)$ of the Orlik-Solomon algebra $A$ of an
arrangement provided with the differential of multiplication by an element
$a\in A$ of degree 1. The second kind is the cohomology $H^*(M,\Lo(a))$
of the complement $M$ of the arrangement with the coefficients in
the local system $\Lo(a)$ defined
by the 1-differential form corresponding to $a$.
According to \cite{ESV,STV} these cohomologies coincide if $a$ is sufficiently
generic, more precisely if the residue of the differential form
at certain intersections of hyperplanes
is not a positive integer. The problem of vanishing of the cohomology of either
kind also was studied in several papers (e.g. see \cite{Ko, Yu}). 

M.Falk in
\cite{Fa} looked at $H^1(A,a)$ from a different point of view. 
He defined the variety $R_1$ of all the elements
$a\in A_1$ with $H^1(A,a)\not=0$ and 
considered its isomorphism class as an (full) invariant of $A$ 
(its quadratic closure). In \cite{L}, the
characteristic varieties of algebraic curves were used for the case of
arrangements of lines in the projective plane.
For this case, a new sufficient condition was obtained for
the isomorphism of the two kinds of cohomology above. This method was further
developed in \cite{CS}. In particular, the conjecture of Falk
that $R_1$ is the union of linear subspaces of $A_1$
 was proved there.

As for the relation between the two types of the cohomology we
propose the following conjecture (Conjecture \ref{the}).

For every $p$ and almost every $a\in A_1^*$ 
among those with $H^p(M,\Lo(a)) \ne 0$ one has 
$$\dim H^p(M,\Lo(a))=\sup_{N\in\ints^n}\dim H^p(A^*,a+N)$$
where the identification $A_1={\cn}^n$ is given by the natural basis of $A_1$ (in
one-to-one correspondence with the set of hyperplanes). 
Choice of precise meaning of ``almost  all'' can be made stronger (resp. 
weaker), e.g. ``for all $a$ for which the corresponding local system is
outside of a finite (resp. Zariski closed) set in the space of local 
systems with non-vanishing $p$th cohomology''.
The fact that the left
hand side is not smaller than the
 right hand side is a simple corollary of \cite{ESV}
(see Corollary \ref{stronger}).
We prove the stronger version of the conjecture for $p=1$ and $a\in R_1$
in Theorem \ref{proof}.

To address this conjecture it is natural to try to find out more
about the components of $R_1$.
In \cite{Fa}, a basic piece of data for describing a component was so called
neighborly partition of the set of hyperplane. Then the component was given by
a system of linear and quadratic equations. We suggest to start from a linear
system.
Let us associate with every $a\in A_1$ such that 
$\sum_{i=1}^na_i=0$ the set $\X=\X(a)$
of all the elements $X$ of
rank 2 of the intersection lattice of the arrangement such that
$\sum_{H_i\supset X}a_i=0$. Treating elements
 $X$ as sets of hyperplanes intersecting at
$X$ we can describe $\X$ by its incidence matrix $J$.
Our main tool is 
the symmetric matrix $Q=J^tJ-E$ where $E$ has all the entries
1. Since all the off-diagonal entries of $Q$ are either 0 or -1 we can apply the
Vinberg classification and its generalization (see Theorem \ref{vinberg})
to it. Since besides any 1-cocycle
of $A$ is in the null-space of $Q$ we obtain some new
information about $H^1(A,a)$.
The main piece of it is that $\dim H^1(A,a)$ is defined
by $\X(a)$ (see \ref{depend}).
 This implies that the set of irreducible components of $R_1$ is in a one-to-one
correspondence with the set of null-spaces of the matrices
$Q$ ($a\in A_1$) and distinct components intersect only at 0.
In fact, a component $V$ coincides with the set of cocycles for any complex
$(A,a)$ with $a\in V\setminus\{0\}$.
 Moreover, if $a\in R_1$ then the respective matrix
$Q$ should have at least three irreducible components of affine type (see
Section 2 for definition). 
This condition turns out to be quite restrictive and allows us to
classify certain classes of lattices (matroids) with $R_1\not=0$ (see Section
6).

On the other hand, intersecting a projective hyperplane arrangement by 
a plane in general 
position, i.e. using a Lefschetz type argument referred to above, 
 we obtain an arrangement of lines in ${\bf P}^2$. Blowing up ${\bf
P}^2$ at the set $\X'$ of
points corresponding to $\X$ gives a variety $P'$ that
carries quite a lot of information about the initial arrangement.
For instance, we recover $Q$ as the minus
intersection form on $P'$.
To each component of $R_1$ having a positive
dimension we can associate a pencil of curves, i.e. an algebraic map 
of the complement of $\X'$
onto the complement of a certain finite subset of
${\bf P}^1$. The structure of the matrix $Q$ 
is closely related to the structure of this pencil. 
The xistence of the pencils imposes restrictions on 
the euler characteristic of $P'$
which in turn yields strong restrictions on the arrangement.
We obtain the classification and nonexistence results mentioned above 
using a combination of combinatorial and topological points of view.  
 
The layout of the paper is as follows. In section 2, we study properties of
matrices 
$Q$ using the Vinberg classification and generalizing it.
In section 3, we use these properties to
obtain information about $H^1(A,a)$. In section 4, we discuss relations between
$H^p(A,a)$ and $H^p(M,\Lo(a))$ and state the conjecture. In section 5, we
generalize results from \cite{A} and prove a version of the conjecture
for $p=1$. Sections 6 considers from combinatorial and section 7 
from algebro-geometric points of view examples and  classification of some
types of $Q$ both realizable and non-realizable by matroids and arrangements.
In particular we obtain complete classification of arrangements with 
positive dimensional $R_1$ and with each line containing exactly three 
points, the arrangements with matrix $Q$ satisfying certain
properties and bound on the number of lines in an arrangement with 
fixed $r,k$ where $r=\dim R_1$ and  $k$ is the maximal number of
vertices on each line

We thank M.Falk and D.Cohen for a discussion of the first version of our
conjecture.
The second author is grateful to R.Silvotti and H.Terao conversaions with whom
have rekindled his interest in the subject.

\section{Matrices of collections of subsets}
\bigskip
In this section, $I$ is a finite set and $\X$ is a non-empty  collection of its
subsets such that
 each two elements of $\X$ have at most one common element.
We will always assume that both $I$ and $\X$ are linearly ordered.
Let $J=J(\X)$ be the incidence matrix of $\X$,
 i.e. the $|\X|\times |I|$-matrix with 
entries $a_{X,i}=1$ if $i\in X$ and $a_{X,i}=0$ otherwise. Also let $E$ be the 
$|I|\times |I|$-matrix whose every entry is 1. Our goal in this section
 is to analyze the matrix

$$Q=Q(\X)=J^tJ-E.$$

For each real matrix $R$ denote by $V(R)$ the null space of $R$ (consisting of
real column
vectors). Notice that 
$V(E)=\{u\in \rn^{I}\vert \sum_{i\in I}u_i=0\}$. Then for any subspace $V\subset
\rn^I$
put $V^*=V\cap V(E)$. The following observation is immediate.

\begin{lemma}
\label{QJ}
$V(Q)^*=V(J)^*$.
\end{lemma}

Clearly $Q$ is symmetric and its every entry off the diagonal is either -1 or 0.
The decomposition of $Q$
in the direct sum of its indecomposable principle submatrices
defines a partition $\Pi$ of $I$. This
partition 
can be defined in terms of $\X$ as follows. Let $\Gamma$
 be the graph on $I$ whose edges
are the pairs $i,j$ that are not in any $X\in\X$. Then $\Pi$ is the partition of 
$\Gamma$ into its connected components. 

We have $Q=\oplus_{K\in\Pi}Q_K$ where each $Q_K$ is a integer matrix satisfying 
the same conditions as $Q$ and besides being indecomposable. Thus one can apply 
the Vinberg classification (e.g., see \cite{Ka}, pp. 48-49) to $Q_K$.

Let us recall Vinberg's result.
Suppose $R=(a_{ij})$ is a $m\times m$ real indecomposable matrix with two extra
properties: 
(a) $a_{ij}\leq 0$ for $i\not=j$ and (b) $a_{ij}=0$ implies $a_{ji}=0$.
We call a real column vector $u^t=(u_1,\ldots,u_m)$ positive and
write $u>0$ if all $u_i>0$ ($u<0$ if $-u>0$).

Then one and only one of the following three
possibilities holds for $R$:

(i) $R$ is of finite type, i.e., it is positive definite (equivalently for every
$u>0$ we have $Ru>0$);

(ii) $R$ is of affine type, i.e., it is positive semidefinite of rank $m-1$ and its
null space is spanned by a positive vector;

(iii) $R$ is of indefinite type, i.e., there exists a positive
vector $u$ such that $Ru<0$.

Applying that classification to $Q_K$ we will call it simply
finite, affine or indefinite if it is of the respective type.

The following remark will be useful.
Suppose $\X$ does not cover $I$. Then at least one row of $Q$ consists of -1
whence $Q$ is irreducible and indefinite. If $\X$ covers $I$ then all diagonal
elements of $Q$ are non-negative.

Now we can prove the main result of the section.

\begin{proposition}
\label{vinberg}
For matrix $Q(\X)$ there are precisely two possibilities:

(i) For every $K\in\Pi$, $Q_K$ is finite or affine;

(ii) For a unique $K_0\in\Pi$, 
$Q_{K_0}$ is indefinite and $Q_K$ is finite for every other $K\in\Pi$.
\end{proposition}

\proof
It suffices to show that there cannot be distinct $K_1,K_2\in\Pi$ such that 
$Q_1=Q_{K_1}$ is indefinite and $Q_2=Q_{K_2}$ is either indefinite or affine.
Suppose to the contrary
that $Q_i$ ($i=1,2$) are such. Then there exist column vectors
$u^i>0$ in ${\rn}^{K_i}$ such that $Q_1u^1<0$ and  $Q_2u^2<0$ or $Q_2u^2=0$.
 We can also assume that $\sum_{j\in K_1}u^1_j=\sum_{j\in K_2}
u^2_j$. Define $u\in {\rn}^I$ by $u_j=u^1_j$ if $j\in K_1$, $u_j=-u^2_j$
if $j\in K_2$ and $u_j=0$ if $j\in I\setminus(K_1\cup K_2)$. Notice that
$\sum_{j\in I}u_j=0$ whence  $Eu=0$.
Then we have using $(\ ,\ )$ for the standard dot product

$$0\leq (Ju,Ju)=(Qu,u)+(Eu,u)$$
$$=(Q_1u^1,u^1)+(Q_2u^2,u^2)\leq (Q_1u^1,u^1)<0$$
which is a contradiction.                         \qed

\medskip
\begin{definition}
The collection $\X$ is called either 
affine or indefinite if respectively possibility
(i) or (ii) from Proposition \ref{vinberg} realizes for $Q(\X)$.
In the former case, we denote by $\Pi_1(\X)$ the set of elements $K\in\Pi(\X)$
such that $Q_K$ is affine. In the latter case, $K_0\in\Pi(\X)$ is always the
unique element such
that $Q_{K_0}$ is indefinite.
\end{definition}

\medskip
\begin{corollary}
\label{from vinberg}
Put $V=V(\X)=V(Q(\X))$.

(i) If $\X$ is affine then $\dim V=|\Pi_1(\X)|$  and there is a
basis of $V$ consisting of vectors $u^K$ ($K\in \Pi_1(\X))$
 such that the restriction of $u^K$ to $K$ is positive and its restriction 
to the other elements of $\Pi(\X)$ is 0. In particular, $\dim V^*=|\Pi_1(\X)|-1$.

(ii) If $\X$ is indefinite then
 all the vectors from $V$ are 0 on every $K\in \Pi(\X)\setminus\{K_0\}$.
\end{corollary}
In the latter case of Corollary \ref{from vinberg},
 $\Pi(\X)$ does not define $\dim V$.

\noindent We will also use the following observation.

\begin{proposition}
\label{XK}
If $\X$ is affine then for every $X\in\X$ we have $X\cap K\not=\emptyset$ either
for all $K\in\Pi_1(\X)$ or for none of them.
\end{proposition}

\proof
Suppose to the contrary that there exists $X\in\X$ and $K'\in \Pi_1(\X)$ such that
$X\cap K'=\emptyset$ but $X\cap K\not=\emptyset$ for $K\in\bar\Pi\subset\Pi_1(\X)$
with a non-empty $\bar\Pi$. 

First notice that $|\bar\Pi|>1$. Indeed if $\bar\Pi=\{K\}$ then for 
any $u\in\rn^I$
such that $\sum_{i\in X}u_i=0$ its restriction to $K$ would be 0. But there exists
$v\in V(\X)^*$ non-zero on $K$ and $K'$ which contradicts Lemma \ref{QJ}.
Thus there exist distinct $K_1,K_2\in\bar\Pi$.

Now compare $V=V(\X)^*$ and $V'=V(\X')^*$  where $\X'=\X\setminus\{X\}$.
On one hand, since one equation is deleted we have $V\subset V'$. On the other
hand, $\Pi(\X')$ can be obtained from $\Pi(\X)$  
by gluing together the elements intersecting with $X$. Since $K'$ is untouched,
$Q_{K'}$ is unchanged whence
$\X'$ is affine. Since $K_1$ and $K_2$ get glued $|\Pi_1(\X')|<|\Pi_1(\X)|$
whence by Corollary \ref{from vinberg} $\dim V'<\dim V$. The contradiction
completes the proof.                            \qed

\bigskip
\section{$H^1(A,a)$}
\bigskip

In this section, $\A=\{H_1,\ldots,H_n\}$ 
is a linear arrangement of hyperplanes in a space of dimension $\ell$ over an
arbitrary field
and $L$ is its intersection lattice.
 By $L(k)$ we denote the subset of $L$ of all elements of rank 
(i.e., codimension) $k$.
Every $X\in L$ defines the subarrangement $\A_X=\{H\in\A\vert
H\supset X\}$.
We will treat elements $X\in L$ also as subsets of $\bar n=\{1,\ldots,n\}$
putting $i\in X$ if $H_i\in\A_X$. In this sense, we put $L'(2)=\{X\in L(2)\vert
|X|>2\}$. Notice that $L'(2)$ is the set of all $X\in L(2)$ with
indecomposable $\A_X$.
For every two distinct $i,j\in\bar n$ there exists a unique
$X_{ij}\in L(2)$ such that $i,j\in X_{ij}$.

Now let us fix an infinite field $F$ of a characteristic different from 2
and recall the definition of the Orlik-Solomon algebra $A=A(\A)$ (depending on
$L$ and $F$ only). Let $E$ be the graded exterior algebra over $F$
with generators
$e_1,\ldots,e_n$. Define the linear map $\partial:E_p\to E_{p-1}$ by
$$\partial(e_{i_1}\cdots e_{i_p})=\sum_{j=1}^p(-1)^{j-1}e_{i_1}\cdots \hat e_{i_j}
\cdots e_{i_p}.$$
Then  $A$ is the factor algebra
of $E$ by the homogeneous ideal generated by $\partial(
e_{i_1}\cdots e_{i_p})$ for all dependent sets $\{H_{i_1},\ldots,H_{i_p}\}$.
The grading on $E$ induces a grading $A=\oplus_pA_p$ on $A$.
We keep the notation $e_1,\ldots, e_n$ for the images of $e_i$ in $A$. Notice
that these elements form a basis of $A_1$
that allows us to identify $A_1$ with $F^n$.
As in the previous section, we put $A_1^*=\{a\in A_1 
 |\sum_{i=1}^n a_i=0\}$ and $U^*=U\cap A_1^*$ for every subset $U$ of $A_1$.
For every $a\in A_1$ the multiplication by $a$ defines the differential of
degree 1 on $A$. The respective cohomology is denoted by $H^*(A,a)$.

The ultimate goal of this section is to study $H^1(A,a)$ as a function of $a$ and
the set $R_1=\{a\in A_1\vert H^1(A,a)\not=0\}$ (cf. \cite{Fa}). For that we put
  $$Z(a)=\{x\in A_1|ax=0\}$$ 
and notice that $H^1(A,a)=Z(a)/F a$. 

The following three lemmas are spread through the literature
(e.g. see \cite{CS, Fa, L, Yu})
 and straightforward.

\begin{lemma}
\label{rank2}
Let $\rk L=2$ and $a=\sum_{i=1}^na_ie_i\in A_1$ .
 Then $x=\sum_{i=1}^nx_ie_i\in Z(a)$ if and only if
$$a_i\sum_{j=1}^nx_j=(\sum_{j=1}^na_j)x_i,\ {\rm for\ every}\ i=1,\ldots,n.
$$
Therefore

(i) if $a\not=0$ and $\sum_{i=1}^na_i=0$ then $Z(a)$ can be described by the equation
$$\sum_{i=1}^nx_i=0;$$

(ii) if $\sum_{i=1}^na_i\not=0$ then $Z(a)=Fa$;

(iii) if $a=0$ then $Z(a)=A_1$.

If $n=2$ then the conclusion in case (i) coincides with that in case (ii).
\end{lemma}

For every $x=\sum_{i=1}^nx_ie_i\in A_1$ and every $X\in L$ we put
$x(X)=\sum_{i\in X}x_ie_i$.
\begin{lemma}
\label{projection}
Let $\A$ be arbitrary and $a\in A_1$. 
Then $x\in Z(a)$ if and only if $x(X)\in Z(a(X))\subset A_1(\A_X)$ for every
$X\in L(2)$.
\end{lemma}

Combining the two previous lemmas we obtain the following.

\begin{lemma}
\label{system}
Using the notation of Lemma \ref{projection},
$Z(a)$ can be described by the following system:

for every $X\in L'(2)$ such that $a(X)\not=0$ and $\sum_{i\in X}a_i=0$,
$$\sum_{i\in X}x_i=0;\eqno(3.1)$$

for every other $X\in L(2)$ and
every pair $i<j$ from $X$,
$$a_ix_j-a_jx_i=0.\eqno(3.2)$$
\end{lemma}

\medskip
It is known that if $a\in A_1\setminus A_1^*$ then $H^*(A,a)=0$ (\cite{Yu}).
Therefore we restrict our considerations to $a\in A_1^*$. Since the algebra $A$ is
skew commutative $b\in Z(a)$ if and only if $a\in Z(b)$. Thus $Z(a)^*=Z(a)$ for
every $a\in A_1^*\setminus\{0\}$.

Now we introduce the main objects of the section.
For any $a\in A_1^*$ define the subset of $L'(2)$ via
$$\X(a)=\{X\in L'(2)\vert\sum_{i\in X}a_i=0,\ a(X)\not=0\}$$  
and the subset of $\bar n$ via 
$$I(a)=\bigcup_{X\in\X(a)}X\subset\bar n.$$
Notice that if $\A$ is a line arrangement then $\X(a)$ is a set of points and
$I(a)$ is the subarrangement of all the lines passing through any of the points
from $\X$.

The following theorem is the main result of this section.
In it we apply the notation from Section 2 
to the collection (in fact, covering) $\X(a)$ on $I(a)$.
In particular $V(\X(a))$ denotes 
the subspace of $A_1$ of all elements that are 0 outside of $I(a)$ and
whose restrictions to $I(a)$ lie in the null space of $Q=Q(\X(a))$.

\begin{theorem}
\label{main}
Let $a\in A_1^*\setminus\{0\}$.
If $\X(a)=\emptyset$ or $\X(a)$  on $I(a)$  is indefinite
then $a\not\in R_1$.  If this collection is affine then $Z(a)=V(\X(a))^*$,
in particular 
$\dim Z(a)=|\Pi_1(\X(a))|-1$.
\end{theorem}

\proof
Let us keep in mind that for every
 $j\in\bar n\setminus I(a)$ 
 and $i\in\bar n$ we have $X_{ij}\not\in\X(a)$ whence (3.2) holds for every $x\in
Z(a)$.
Thus if there exists
 $j\in\bar n\setminus I(a)$ 
such that $a_j\not=0$ then
$\dim Z(a)=1$ and $a\not\in R_1$.
In particular this covers the case where $\X(a)=\emptyset$.

From now on we can assume that $\X(a)\not=\emptyset$ and
 $a_j=0$ for every $j\not\in I(a)$. For every $j\not\in I(a)$,
$i$ with $a_i\not=0$ and $x\in Z(a)$,
the equalities (3.2) imply
 $x_j=0$. Thus $x$ can be non- zero  
only on $I(a)$ and we can identify $x$ with its restriction to 
 $I(a)$. 
By definition of $\X(a)$ and Lemma \ref{system}, (3.1) holds
for every $X\in\X(a)$ whence $J(\X(a))x=0$. 
Since besides 
$\sum_{i=1}^nx_i=0$ we have $Q(\X(a))x=0$, i.e., $x\in V(\X(a))^*$.
  Thus $Z(a)\subset V(\X(a))^*$.

Now suppose that $\X$ is indefinite
 and $x\in Z(a)$. Since $x\in V(\X(a))$ we conclude that $x$ is 0 on each
$K\in\Pi=\Pi(\X(a))$ such that $K\not=K_0$. Let us study the coordinates 
of $x$ on $K_0$. If
$i,j\in K_0$ are such that $X_{ij}\not\in\X(a)$ then (3.2) holds. Since $K_0$ is
an element of $\Pi$  and $a\not=0$, (3.2) holds for every $i,j\in K_0$ and 
$x=ca$ 
on $K_0$ for some $c\in F$. Summing up we see that
$x=ca$ whence $\dim Z(a)=1$, 
i.e., $a\not\in R_1$.

Now suppose $\X=\X(a)$ is affine and $x\in V(\X)^*$, i.e., $Q(\X)x=0$
and $Ex=0$. By Lemma \ref{QJ} we have $J(\X)x=0$,
i.e.,  (3.1) holds for $X\in\X$.
Moreover for every $i,j\in I(a)$ such that $X_{ij} \not\in\X$
 we have 
$i,j\in K$ for some $K\in\Pi$. If $Q_K$ is finite then $x_i=x_j=0$ and (3.2)
holds tautologically. If $Q_K$ is affine then the kernel of $Q_K$ 
is one-dimensional and (3.2) holds again. 
By Lemma \ref{system} $x\in Z(a)$
that completes the proof.                      \qed

\medskip

\begin{corollary}
\label{depend}
For every $a\in A_1^*$ the dimension of
 $Z(a)$ is defined by $\X(a)$.
If the collection $\X(a)$ is affine (in particular if $a\in R_1$)
then the space $Z(a)$ itself is defined by $\X(a)$. If besides $a\not=0$ then for
every $b\in Z(a)\setminus\{0\}$ we have $Z(b)=Z(a)$.
\end{corollary}
\proof
Only the last statement needs a proof. Let $b\in Z(a)=V(\X(a))^*$ and $b\not=0$.
Then there exists $K\in\Pi_1(\X(a))$ such that $b_i\not=0$ for every $i\in K$.
Take $X\in\X(a)$. We have $\sum_{i\in X}b_i=0$. Besides by Proposition \ref{XK},
since $a\ne 0$ we have
$X\cap K\not=\emptyset$ whence $b(X)\not=0$. This implies that $X\in\X(b)$ and
$\X(a)\subset\X(b)$. By symmetry $\X(a)=\X(b)$ and the result follows.
             \qed 

\medskip
\begin{corollary} (cf. \cite{L,CS})
\label{components}
All the  irreducible components of $R_1$ 
are linear.
\end{corollary}
\proof
Due to Theorem \ref{main},  $R_1$ is the union 
of finite number of linear spaces. The result follows.       \qed

\medskip
The previous corollary leaves open the question
 what the linear components of $R_1$ are. The answer turns out to be simple.

\begin{corollary}
\label{maximality}
For every $a\in R_1\setminus\{0\}$, the space $Z(a)=V(\X(a))^*$
 is a maximal linear space in $R_1$ whence an irreducible component of $R_1$.
\end{corollary}
\proof
Theorem \ref{main} and Corollary \ref{depend} imply that
$$R_1=\bigcup_{a\in R_1\setminus\{0\}}Z(a)$$ 
and each two distinct $Z(a)$ and $Z(b)$
intersect only at 0. Since there is only finite number of these spaces and
$|F|=\infty$,
each one of them is a maximal
linear subspace in $R_1$.
                     \qed

\medskip
\begin{corollary}
\label{intersection}
The pairwise intersection of the irreducible 
components of $R_1$ is 0.
\end{corollary}

\medskip
Corollary \ref{maximality} gives a description of the irreducible components of
$R_1$ in terms of elements of $A_1^*$. One can give a description by properties
of subsets of $L'(2)$. For every $\X\subset L'(2)$ put $I(\X)=\cup_{X\in\X}X$
and as above denote by $V(\X)^*$ the subspace of $A_1^*$ corresponding to $\X$ on
$I(\X)$.

\begin{theorem}
\label{irreducible}
Let $\X\subset L'(2)$ satisfy the following conditions. (i) It is
 affine on $I=I(\X)$, (ii) $|\Pi_1(\X)|\geq 3$ and 
(iii) no $X\in\X$ lies in $\bigcup_{K\in\Pi\setminus\Pi_1}K$. 
Then $V(\X)$ is an irreducible component of $R_1$. The map $\X\mapsto
V(\X)$ 
defines a one-to-one
correspondence of the set of collections satisfying conditions (i)-(iii)
with the set of irreducible components of $R_1$.
\end{theorem}

\proof
Suppose $\X$ satisfies the conditions (i)-(iii). Put $V=V(\X)^*$ and
fix $a\in V\setminus\{0\}$. If $X\in\X$ then  $\sum_{i\in X}a_i=0$. Besides
since $a\not=0$ and $\X$ satisfies (iii), Proposition \ref{XK} implies that
$a(X)\not=0$. Thus $X\in\X(a)$ and $\X\subset\X(a)$. 

On the other hand, let $X\in \X(a)$. Since $a(X)\not=0$ there exists
$K\in\Pi_1(\X)$ such that $X\cap K\not=\emptyset$. Since $\sum_{i\in X\cap I}a_i=
0$ and $a(K)>0$ or $a(K)<0$ there exist distinct $K_1,K_2\in\Pi_1(\X)$ 
intersecting with $X$.
 This means that $X=X_{ij}$ for some $i\in K_1$ and $j\in K_2$.
Since by definition of $\Pi(\X)$ we have $X_{ij}\in\X$ it follows that $X\in\X$
whence $\X(a)=\X$. 

Now Theorem \ref{main} imply that $V=V(\X(a))^*=Z(a)$.
By condition (ii) $\dim Z(a)\geq 2$ whence $a\in R_1$.
 Then Corollary \ref{maximality}
implies that $V(\X)$ is an irreducible component of $R_1$.

Conversely let $V$ be an irreducible component of $R_1$ and $a\in
V\setminus\{0\}$. By Corollaries \ref{maximality} and \ref{intersection},
$V=V(\X(a))^*$. Clearly $\X(a)$ satisfies the conditions (i)-(iii) above. Since
$\X(a)\mapsto V$, the map is surjective. 

Finally suppose that $V=V(\X_1)=V(\X_2)$ where $X_i$ ($i=1,2$) satisfy the
conditions (i)-(iii). Choose a non-zero $a\in V$. By the first part of the proof,
 $\X_1=\X(a)=\X_2$ whence the above map is injective. This completes the proof.
                    \qed

\medskip
\begin{remark}
\label{falk}
Let us compare the description of the irreducible components of $R_1$ from
Theorem \ref{irreducible} and that from \cite{Fa}. Fix an $\X$ as in the above
theorem and consider the submatroid of $L$ on $I'=\bigcup_{K\in\Pi_1(\X)}K$.
Then it is easy to see that the partition $\pi$
generated by $\Pi_1(\X)$ on $I'$ is
neighborly (using the terminology from \cite{Fa}, p.144)
 and the set of all polychrome
flats of rank 2 coincides with $\X$ (moreover every polychrome $X$ intersects
every element of $\pi$). Thus $L_{\pi}$ from \cite{Fa} coincides
with $V(\X)$ which implies $L_{\pi}=V_{\pi}$ (cf. \cite{Fa}, Remark 3.15). 
\end{remark}

\bigskip

\section{$H^*(A^*,a)$ and cohomology of local systems}
\bigskip
For the rest of the paper we consider only complex arrangements and $F=\cn$.
Let $\A=\{H_1,\ldots,H_n\}$
be an arrangement in $\cn^{\ell}$ and $A$ its Orlik-Solomon algebra.
According to the Orlik-Solomon theorem \cite{OS},
$A$ is isomorphic to 
the algebra of closed differential forms on the complement
 $M'$ of $\A$ generated by the forms
$\omega_i={d\alpha_i\over{\alpha_i}}$ 
where $\alpha_i$ is any linear functional
on $\cn^{\ell}$ with kernel $H_i$. 
The isomorphism is given by
 $e_i\mapsto \omega_i$.
Due to the Brieskorn theorem \cite{Br} this
algebra of forms is isomorphic to $H^*(M',\cn)$ under the de Rham homomorphism.
If we projectivize $\A$ and denote by $M$ its complement in
${\bf P}\cn^{{\ell}-1}$ then $H^*(M,\cn)$ is isomorphic to the subalgebra $A^*$ of
$A$ generated by $A_1^*$. 

Denote the form on $M$ corresponding 
 to $a\in A_1^*$ by $\omega(a)$. Since $\omega(a)$ is closed and $\omega(a)\wedge
\omega(a)=0$ it
 defines a local system on $M$ that we denote by $\Lo(a)$.
More explicitly $\Lo(a)$ is associated with the one-dimensional representation of
$\pi_1(M)$ sending its generator corresponding to $H_k$ to exp$(2\pi ia_k)$.

 We will need the following theorem.

\medskip
{\bf Theorem (STV)} \cite{STV}. 
{\it Suppose that $a\in A_1^*$ and
for all $X\in L$ such that $\A_X$ is indecomposable, the sum
 $\sum_{i\in X}a_i$ 
 is not a positive integer. Then 
$$H^p(M,\Lo(a))=H^p(A^*,a)$$
for every $p$.}

\medskip
Let us record a simple but important observation for the rest of 
this section.

\begin{lemma}
\label{shifts}
For every $a\in A_1^*$

(i) $H^*(A^*,a)=H^*(A,\la a)$ for every  $\la\in\cn^*$;

(ii) $H^*(M,\Lo(a))=H^*(M,\Lo(a+N))$
for every $N=(N_i)\in\ints^n$ such that 
$\sum_{i=1}^nN_i=0$.
\end{lemma}

\medskip
The main result of this section is as follows.
\begin{proposition}
\label{inequality}
For every $a\in A_1^*$ and every $p$
$$\dim H^p(M,\Lo(a))\geq \dim H^p(A^*,a).$$
\end{proposition}

\proof
Let $a\in A^*_1$.
Fix a real positive number $\epsilon$ such that for 
every $X\in L$ and $0< t<\epsilon$ the number $(1-t)\sum_{i\in X}a_i$ 
is not a positive integer.
Then by Theorem (STV) and Lemma \ref{shifts}(i) we have for every $p$
$$H^p(M,\Lo((1-t)a)= H^p(A^*,(1-t)a)=H^p(A^*,a).$$
Using the upper semicontinuity of the dimension of cohomlogy for $t\to 0$
we obtain the result.              \qed

\medskip
\begin{corollary}
\label{stronger}
For every $a\in A_1^*$, $N\in\ints^n$ and $p$

$$\dim H^p(M,\Lo(a))\geq\dim H^p(A^*,a+N).$$
\end{corollary}

\medskip
Corollary \ref{stronger} and examples justify the following conjecture.

\begin{conjecture}
\label{the}
For every $p$ and almost all $a\in A_1^*$ among those with \hfill\break
$H^p(M,\Lo(a))\ne 0$ one has
$$\dim H^p(M,\Lo(a))=\sup_{N\in\ints^n}\dim H^p(A^*,a+N).$$
\end{conjecture}

\medskip
See Introduction for a discussion of ``almost all''.
A version of this conjecture will be proved for $p=1$ in the next
section.
\bigskip
\section{Characteristic varieties}
\bigskip

Our first goal in this section is to give
a more precise statement of the Proposition
1.7 from \cite {A}. Since we work here in more generality than in 
the previous sections we start by fixing notation.

Let $M$ be an arbitrary quasiprojective variety. We identify the torus 
$H^1(M,{\cn}^*)=\Hom (\pi_1(M),{\cn}^*)$ with the set of all (rank one)
local systems on $M$. We also call $H^1(M,\cn^*)$ the torus of characters 
$\Char(\pi_1(M))$ of $\pi_1(M)$. Let $\Sigma^k(M) \subset\Char(\pi_1(M))$ 
be the  subset consisting
of those local systems $\Lo$ for which $\dim H^1(M,\Lo)
\ge k$. Then $\Sigma^k(M)$ is an algebraic subvariety of $\Char(\pi_1(M))$. 
We want to study the structure of irreducible components of $\Sigma^k(M)$ of
positive dimension. According to \cite{A},
such a component belongs to a component of 
$\Sigma^1(M)$. The latter  defines a map 
$f: M \rightarrow C$ to a smooth curve $C$ and has a form
 $\rho f^*(H^1(C,{\cn}^*))$ for some torsion point $\rho$ in
 $\Char(\pi_1(M))$.
 We assume that $C$ {\it isn't} compact (which is always the case if 
$M$ is a complement to an arrangement, cf. 5.2),
i.e. $\pi_1(C)$ is a free group with, say, $n$ generators. We shall 
denote such a group below as $F_n$.
Notice that for a local system $\Lo$ on $C$ we have
$\dim H^1(C,\Lo)=n-1$ if $\Lo$ is not trivial and
for the trivial local system it is $n$ (indeed $F_n$ has cohomological dimension 
equal to $1$, the euler characteristic is independent of a local system 
and is equal to $1-n$ and $\dim H^0(C,\Lo)$ is 1 or zero depending on whether the
system is trivial or not). 

\begin{theorem} 
\label{arapura}
Let $V$ be a component of $\Sigma^k(M)$ containing the identity
character. Then there 
is a finite set of local systems $E \subset V$ (possibly empty) and an
integer $s$ such that for any local system $\Lo$ in $V\setminus E$
we have $\dim H^1(M,\Lo)=s$ (i.e. $max \{ k \vert V \subset \Sigma^k(M)\}=s$).
For such a component and the corresponding 
map $f: M \rightarrow C$ such that $V=f^*(H^1(C,{\cn}^*))$ the
number of generators of the fundamental group of $C$ is $s+1$.
\end{theorem}

\proof 
By theorem 1.6 in \cite{A} any local system in $V$ 
 has form $f^*(\Lo)$ for an appropriate local system on $C$. 
The main step in the argument is the following statement 
implicitly, as D.Arapura pointed out, contained in \cite{A}. For all
but finitely many local systems in $V$ we have 
$$H^1(M,f^*(\Lo))=H^1(C,\Lo). \eqno (5.1)$$
Using that
$C$ has cohomological dimension
1, the projection formula \hfill\break
$R^qf_*(f^*(\Lo))=R^qf_*{\cn} \otimes\Lo$ and
the Leray spectral sequence 
$$E_2^{p,q}=H^p(C,R^qf_*(f^*(\Lo)))\Longrightarrow H^{p+q}(M,f^*(\Lo))$$ 
we see
that (5.1) is a corollary of the vanishing of 
$H^0(C,R^1f_*({\cn})\otimes\Lo)$ for all but finitely many local systems
$\Lo$ on $C$. As it is shown in \cite{A}, if $U$ is the subset of $C$ 
over which $f$ is a locally trivial fibration, then 
 $$H^0(C,R^qf_*{\cn} \otimes\Lo)=H^0(U,R^qf_*{\cn} \otimes\Lo).$$
let $g_i$ ($i=1,...,s$) be a system of generators of $\pi_1(U)$ and  
$\lambda_i, i=1,...,s$ (resp. $\mu_{i,j}, i=1,..,s,
j=1,\ldots,rkR^1f_*{\cn})$) the eigenvalues of the action of $g_i$ on $\Lo$ 
(resp. on $R^1f_*({\cn})$). Then $H^0(U,R^qf_*{\cn} \otimes\Lo)=0$
unless for each $i$ there is $j$ such that 
$\lambda_i \cdot \mu_{i,j}=1$. 
Hence (5.1) is valid for all but at most 
$\dim (R^qf_*{\cn}) \cdot s$ local systems.                  \qed

\medskip
Now let $\A$ be an arrangement of $n$ lines
 in ${\bf P}\cn^2$, $M={\bf P}\cn^2\setminus\bigcup_{H\in\A}H$
and $A^*=H^*(M,\cn)$ the subalgebra of the Orlik-Solomon algebra generated by
$A_1^*$. 
We have the following extension of results of \cite{L}
to the case of general arrangements (not necessarily containing a line
transversal to the rest of the lines of the arrangement) (cf. also \cite{CS}).

\begin{theorem}
\label{cover}
Let $V^k\subset A_1^*$ be an
irreducible component of  
the variety of those $a\in A_1^*$ for which 
$\dim H^1(A^*,a) \ge k$. 
 Then the image of ${V^k}$ 
under the map ${\rm exp}:a\mapsto \Lo(a)$
 is an irreducible component of 
$\Sigma^s (M)$  and exp
is a universal covering of this component.
\end{theorem}

\proof 
By Theorem \ref{main} and Corollary \ref{maximality} $V^k$ is a linear space.
By Proposition \ref{inequality} we have $\dim H^1(A^*,a) \le \dim H^1(M,\Lo(a))$.
Thus ${\rm exp}V^k \subset\V$  where $\V$ is an irreducible 
component of ${\Sigma^s}$ with $s \ge k$. In this case $\V$ is a subtorus of 
$H^1(M,\cn^*)$ (i.e. no translation is needed) since $V^k$ contains the
origin.
Let $f: M \rightarrow C$ be a map as in Theorem \ref{arapura}.
Then $C$ is a rational curve,  i.e., the complement in ${\bf P}^1$ to a 
set $\{p_1,\ldots,p_{s+1}\}$. Indeed, according to H.Hironaka, 
$f$ is a restriction of $\bar f: \bar M \rightarrow \bar C$ where
$\bar M$ and $\bar C$ are compactifications of $M$ and $C$ 
respectively. If $\bar C$ had a positive genus it would admit
a non trivial holomorphic 1-form and its pull back 
via $\bar f$ to $\bar M$ would yield a non trivial holomorphic 1-form on
$\bar M$ contradicting the rationality of $\bar M$.

Denote by $\omega_i$ a 1-form on ${\bf P}^1$ that generates
$H^1({\bf P}^1\setminus\{p_i\},\cn)$ ($i=1,\ldots,s+1$).
The pull backs of these forms 
span an $s$-dimensional space of $Z(a)$.
 Since exp
is a local homeomorphism and $\dim V^k \le s$ the result follows.
              \qed

\medskip
Finally we can prove a version of Conjecture \ref{the} for $p=1$.
Recall that for every $a\in A_1^*$ we have $H^1(A^*,a)=H^1(A,a)$.

\begin{theorem}
\label{proof}
For all but finitely many cosets $mod\  {\bf Z}^n$ of $a\in R_1$ 
$$\dim H^1(M,\Lo(a))=\sup_{N\in\ints^n}\dim H^1(A,a+N).$$
\end{theorem}

\proof 
The arguments are similar to the proof of Theorem \ref{cover}.
We first notice that since
$\dim H^1(A,a)>0$, the local system $\Lo(a)$ 
belongs to a component of the characteristic variety of positive
dimension containing the identity character.
Indeed, since $\dim Z(a)>1$, one can find a deformation $a(t)$ of $a$
($a(0)=a$) for $0\leq t\leq 1$ such that
 $a(t)$ is not proportional to $a(t')$ for $t\not=t'$. Then $a(t)\in R_1$ and
$\Lo(a(t))$ provide a deformation of $\Lo(a)$.

Next let us consider the map $\pi:M\to
{\bf P}^1\setminus \{p_1,..,p_m\}$
 corresponding to a component of the characteristic variety
containing $\Lo(a)$.
Then $\Lo(a)$ is a pullback of a local system 
$\Lo(\omega)$ on ${\bf P}^1\setminus\{p_1,..,p_m\}$ 
corresponding to a differential form $\omega$.
 Let $\tilde \omega$
be the form with constant coefficients on $M$
representing the cohomology class of the pullback of $\omega$
(such a form exists according to Brieskorn's theorem \cite{Br}).
Then $\tilde \omega$ and $\omega(a)$ represent the same local system
and hence $\omega(a)=\tilde \omega+\omega(N)$ for some $N\in{\ints^n}$. Equivalenlty
$a=\tilde a+N$ where $\tilde a\in A_1^*$ and $\tilde\omega=\omega(\tilde a)$.

Now we claim that 
$H^1(A,\tilde a)=H^1(M,\Lo(a))$. 
We need to show that $\dim H^1(A,\tilde a) \ge m-1$
since $\dim H^1(M,\Lo(a))=m-1$ with only finitely many exceptions
by \cite{A}. Indeed the map $H^1({\bf P}^1\setminus\{p_1,..,p_m\})
\rightarrow A_1^*$ is injective (by Leray Spectral
sequence) and its image lies in $Z(\tilde a)$.
This gives us 
an $(m-1)$- dimensional subspace in $H^1(A,\tilde a)$ that completes the
proof.                               \qed.

\medskip Exceptional cosets, for which the equality of the theorem
fails may indeed exist. It follows from Example 4.5  in \cite{CS}.
The arrangement in this example consists of seven lines that can be represented
in three different ways as the union of 
braid arrangement with a 
line connecting two double points. For this arrangement 
$\Sigma^1$ contains non empty $\Sigma^2$ consisting of one point outside
of the origin.

\medskip
Finally let us note that theorems 5.2 and 5.3 
are valid for arbitrary arrangements in ${\bf P}({\cn}^n)$. 
Indeed by Lefschetz theorem (cf. \cite{Mi}) for a generic plane 
$H$ the cohomology of the complement to 
arrangement in $P^n$ and its intersection with H are the same in
dimension 1 and injects in dimension 2. Hence the cohomology   
of the Orlik-Solomon algebra in dimension 1 is the same for 
both arrangements.
On the other hand characteristic varieties also 
are the same since they depend only on the  fundamental group 
(cf. definition in section 1 of \cite{L}). This yields 5.2 for
arrangements in ${\bf P}^n$.
The theorem 5.3 follows from this comparison and the 
Lefschetz theorem for cohomology with coefficients in 
local system (which is a consequence of the same argument as 
in the case of constant coefficients, i.e. that the complement
in ${\bf P}^n$ obtained from the complement in $H$ by attaching cells 
of dimension greater than 2 cf. \cite{Mi},\cite{Si})

\bigskip

\section{Types of affine matrices, realizations and examples}
\bigskip
In this section, we exhibit examples of types of matrices $Q$
 with affine components
and discuss the problem of their realizations by matroids and arrangements.

Recall that for every subset $\X\subset L'(2)$ we put
$I=I(\X)=\bigcup_{X\in\X}X$ and $n=|I|$. Then $Q=Q(\X)=(q_{ij})$ is a symmetric
$n\times n$-matrix with $m_i=q_{ii}\in \ints_+$
and $q_{ij}=-1$ or 0 for $i\not=j$.
To parametrize 
matrices $Q$ we write $Q=Q(m,\Gamma)$ where $\Gamma$ is the graph on
$n$ vertices $v_i,\ldots,v_n$ with edges $\{v_i,v_j\}$ for $q_{ij}=-1$. The vector
$m=(m_1,\ldots,m_n)\in\ints^n_+$ can be regarded as a
labeling of the vertices of $\Gamma$. The decomposition
of $Q$ into the direct sum of indecomposable
components corresponds to the partition
of $\Gamma$ into connected components.

The section is organized in the following way. First we collect some
information and series of examples of pairs $(\Gamma,m)$ with $Q(m,\Gamma)$
affine indecomposable. Then 
we attend the question which of the affine matrices and more importantly
matrices $Q$ having at least three affine components can be realized by
matroids.  These (especially the latter) are more restrictive conditions which
allows us to classify some classes of examples. The final question is whether
the matroids appearing in
this way can be represented by arrangements of complex lines. We give some
examples of representations leaving negative results
till the next section where we will apply different methods.

\medskip
\noindent{\bf 6.1. Affine indecomposable matrices.}

\medskip
\noindent We start with the following two general observations. 

\begin{proposition}
\label{finite}
Let $\Gamma$ be an arbitrary (finite) connected graph. Then
the set $\M(\Gamma)$ of vectors $m\in\ints_+^n$, for which $Q(m,\Gamma)$
 is affine, is finite.
\end{proposition}

\proof
Suppose that $\Gamma$ has $n$ vertices. For each $m\in\M(\Gamma)$ 
denote by $Q_k$ the principal submatrix of $Q=Q(m,\Gamma)$ such that the involved
rows have $k$ smallest $m_i$ and by $D_k$ the determinant of $Q_k$.
 Without loss of generality assume that
$m_1\leq m_2\leq\cdots\leq m_n$. 
\cite{Ka}, Lemma 4.4 implies that for $k<n$ matrix
$Q_k$ is a direct sum of finite matrices whence $D_k>0$. Thus for fixed $k<n$ and
$Q_k$ determinant $D_n$ can be viewed as
a polynomial in $m_{k+1},\ldots,m_n$ with a positive coefficient $D_k$ of the
highest degree term $m_{k+1}\cdots m_n$. Since $D_n=0$ the possible value of
$m_{k+1}$ is bounded from above whence $m_{k+1}$ can afford only finite set of
values. The result follows by induction on $k$.
               \qed

\medskip
\begin{proposition}
\label{parameter}
Let $\Gamma$ be a connected graph on $n$ vertices with a set $E$ of edges. Let
$\N$ be the set of all positive integer column
vectors $u$ ( $u^t=(u_1,\ldots,u_n)$) such that 
$u_i$ are mutually relatively prime and each $u_i$ divides
$u^{(i)}=\sum_{\{i,j\}\in E}u_j$. Define $m=f(u)$ via $m_i=u^{(i)}/u_i$.
Then (i) $m\in\M(\Gamma)$, (ii) $u$ lies in the null space of $Q(m,\Gamma)$ and (iii)
$f:\N\to\M(\Gamma)$ is  a
bijection.
\end{proposition}

\proof
Fix $u\in\N$ and consider the respective matrix $Q=Q(m,\Gamma)$.
Clearly $Qu=0$ whence $Q$ is affine (see \cite{Ka}, Corollary 4.3) which proves
(i) and (ii). Now let $Q$ be an  arbitrary
integer affine matrix. Then its null space contains  a unique positive
integer vector with mutually relatively prime coordinates which implies (iii).
                        \qed
\medskip

\noindent {\bf Examples of affine indecomposable matrices.}

\smallskip
1. {\it
Laplace matrices of graphs.} Here $G$ is an arbitrary connected graph and $Q$ is
its Laplace matrix, i.e., $Q=Q(d,\Gamma)$ where $d_i$ is the degree of $v_i$.
It is easy to check that $Q$ is affine, in particular a column vector $u$
generating its null-space can be taken as $u_i=1$ for all $i$.

2. {\it Affine generalized Cartan matrices.} Here $\Gamma$ is an extended Dynkin
diagram with single edges and $m_i=2$ for all $i$. Notice that for the diagrams of
type other than $A_{\ell}^{(1)}$ the labeling is different from that in the
previous example.

3. {\it Full graphs.} Let $\Gamma$ be the full graph on $n$ vertices.
Then it is not hard to describe all the labeling $m$ so that 
$Q=Q(m,\Gamma)$ is affine.

\begin{proposition}
\label{full graph}
For the full graph $\Gamma$ the matrix $Q=Q(m,\Gamma)$ is affine if and only if 
$$\sum_{i=1}^n{1\over{m_i+1}}=1.\eqno(6.1)$$
\end{proposition}

\proof
Since the null space of $Q$ should contain a positive vector the following system
should have a positive solution:
$$(m_i+1)x_i=\sum_{j=1}^nx_j\ ,1\leq i\leq n. \eqno(6.2)$$
 Clearly (6.1) is equivalent to it.
 Conversely (6.1) implies that any solution of (6.2) is proportional to
$$x_i={1\over{m_i+1}},\ i=1,\ldots,n.$$
                 \qed

\medskip
Notice that the Laplace matrix of the full graph has $m_i=n-1$ for all $i$.

4. {\it Bushes.} Let $\Gamma$ be the bush on $n+1$ vertices. By that, we mean
that $\Gamma$ is a tree
whose one vertex $v_0$ has degree $n$ (the root) and all the other vertices
$v_1,\ldots, v_n$ (changing
enumeration slightly) are leaves.
An argument similar to that of Proposition \ref{full graph}
shows that $Q$ is affine if and only if 
$$m_0=\sum_{i=1}^n{1\over{m_i}}.\eqno(6.3)$$
If (6.3) holds then the null space of $Q$ is generated by the vector
$u^t=(1,1/m_1,\ldots,1/m_n)$.
As a concrete example, the complete list of admissible vectors $m$ for
the extended Dynkin diagram $D_4^{(1)}$
consists of the following 20 vectors (up to permutations of leaves):
(4,1,1,1,1) (the Laplace matrix),
 (3,2,2,1,1), (2,2,2,2,2) (the Cartan matrix), (2,3,3,3,1),
(2,6,3,2,1), (2,4,4,2,1), (1,4,4,4,4), (1,3,4,4,6), (1,12,4,3,3), (1,6,6,3,3),
(1,12,12,3,2), (1,15,10,3,2), (1,18,9,3,2), (1,24,8,3,2), (1,42,7,3,2),
(1,8,8,4,2), (1,12,6,4,2), (1,20,5,4,2), (1,10,5,5,2), (1,6,6,6,2).

\medskip
\noindent{\bf 6.2. Realization of matrices $Q$ by matroids.}

\medskip
\noindent Here we study the problem
 about a realization of $Q$ by a matroid of rank at
most three. A matrix $Q$ is indecomposable
affine in the first part of the subsection and then
will have several affine components.
We describe matroids as in section 3
by the incidence
matrices $J$ of the collections of dependent sets. More precisely the columns
of $J$ are parametrized by the points of the matroid and the rows by its dependent
sets. 
If the matroid is representable over a field $\cn$ by an arrangement of
lines in the projective plane then the columns of $J$ correspond to 
the lines and its rows to the multiple points of their intersections.
Abstractly $J$ is a 0-1 matrix
whose any two columns both have
entries 1 in at most one common row.
The decomposition of $Q$ into indecomposables defines a partition $\Pi$ of
columns of $J$ hence for each $K\in\Pi$ defines a matrix $J_K$ of the
respective 
columns. We will always assume that $J_K$ does not have 0-rows for every $K$.
A realization of $Q$ by $J$ is the identity 
$Q=J^tJ-E$. 

 Not all affine matrices are realizable.

\begin{example}
\label{realize}
Let $Q$ be the Laplace matrix of the graph $\Gamma$ on vertices \hfill\break
$\{v_1,\ldots,
v_5\}$ which is the full graph on $\{v_2,\ldots,v_5\}$ extended by the single edge
$\{1,5\}$. We want to observe that $Q$ is not realizable. Moreover consider
the indecomposable component $S$ (on rows 1 to 4) of $Q+E$:
$$S=\left(\matrix{2&1&1&1\cr
                  1&4&0&0\cr
                  1&0&4&0\cr
                  1&0&0&4\cr}\right).$$
Suppose $S=J^tJ$ for a 0-1 matrix $J$. Then the first column of $J$ has precisely
two 1's in some rows and 
each of the following three columns has 1 in one of these rows.
 Thus some two of these three columns have 1 in
 common row
which is a contradiction. Therefore $S$ and $Q$ are not
realizable.
\end{example}

If an affine matrix is realizable then it could have several different
(up to permutations of rows in $J$) realizations. Clearly the number of
realizations is finite. For instance, the Cartan matrix of type
$D_4^{(1)}$ has 3 different realizations. For certain ordering on the vertices
they are
$$J_i\oplus\left(\matrix{1\cr 1\cr 1\cr}\right)$$
where
$$J_1=\left(\matrix{1&1&1&1\cr
                    1&0&0&0\cr
                    1&0&0&0\cr
                    0&1&0&0\cr
                    0&1&0&0\cr
                    0&0&1&0\cr
                    0&0&1&0\cr
                    0&0&0&1\cr
                    0&0&0&1\cr}\right),
  J_2=\left(\matrix{1&1&1&0\cr
                    1&0&0&1\cr
                    1&0&0&0\cr
                    0&1&0&1\cr
                    0&1&0&0\cr
                    0&0&1&1\cr
                    0&0&1&0\cr}\right),
  J_3=\left(\matrix{1&1&0&0\cr
                    1&0&1&0\cr
                    1&0&0&1\cr
                    0&1&1&0\cr
                    0&1&0&1\cr
                    0&0&1&1\cr}\right).$$

\medskip
More interesting problem for us is the realization of decomposable matrices with
at least 3 affine components.               
This class turns out to be more restrictive. 
While studying this class we sometimes
apply the dual point of view to a matrix $J$ realizing $Q$.
Each column $C_i$ of it can be considered as a subset
of the set $R$ of all rows of $J$.  Since there are no 0-rows in $J_K$,
the system $\C(K)=
\{C_i\}_{i\in K}$ is a covering of $R$ for
each $K\in\Pi$. If $K'\in\Pi$ and $K'\not=K$ then $C_j$ for every $j\in K'$ is a 
{\it transversal} of $\C(K)$. By that we mean that $C_j\cap C_i\not=\emptyset$ for
all $i\in K$ whence $|C_i\cap C_j|=1$.
This immediately implies the inequalities

$$|\bar K|\leq |C_j|\leq |K|\eqno(6.4)$$
where $j\in K'$ and $\bar K$ is a subset of $K$ such that $C_i\cap C_k=\emptyset$
for every $i,k\in\bar K$.

The inequalities (6.4) allow us to classify in some sense
 a series of realizable matrices
with at least 3 affine components.

\begin{theorem}
\label{full graphs} 
Suppose $Q=Q(m,\Gamma)$ is a realizable matrix whose graph
$\Gamma$ has at least two connected components that are full graphs with
at least two vertices each and the matrices corresponding to the components
are affine.
Then there exists an integer $n\geq 2$ such that
all the diagonal elements of $Q$ are $n-1$,
 $|\Pi|\leq n+1$ and $|K|=n$ for each $K\in\Pi$.
If $|\Pi|=n+1$ (i.e. maximal) then
the realizations of $Q$ are parametrized by latin $n\times n$ squares with
entries $1,\ldots,n$.
\end{theorem}

\proof
First notice that for any matrix corresponding to a full graph 
all pairs $\{C_i,C_j\}$ in any realization $J$
are disjoint. Then (6.4) implies that $|C_i|=|K|$ for each $i$ and each
$K\in\Pi$. Put $n=|K|$. Since each $Q_K$ is at least 2$\times
$2-matrix we have $n\geq 2$. Clearly all the diagonal elements of $Q$ are $n-1$.

Now order all the elements of $\Pi$ as $K_0,K_1,\ldots,K_{\ell-1}$. Ordering 
the elements inside $K_s$ ($s=0,1$) and the rows of $J$
 appropriately we can assume without loss of
generality that (i) $C_i=\{n(i-1)+1,\ldots,ni\}$ for $i\in K_0$, i.e., $1\leq i\leq
n$, (ii) $C_j=
\{j-n,j,\ldots,j+(n-2)n\}$ for $j\in K_1$, i.e., $n+1\leq j\leq 2n$, and
(iii) for every $r>0$ the first $n$ rows of $J_r=J_{K_r}$ form the identity matrix.
Now break all the rows of $J$ into blocks of $n$ subsequent rows from $jn+1$ to
$(j+1)n$ ($j=0,1,\ldots,n-1$). All the blocks of $J_r$ ($r=1,\ldots,\ell-1$)
are in bijection with the set of $n$-permutations $\{\sigma(i,j)\}$ ($1\leq
i\leq n,\ 1\leq j\leq\ell-1$) with the following properties:

(i) $\sigma(1,j)=\sigma(i,1)=e$ (the identity permutation) for all $i,j$;

(ii) every two permutations $\sigma(i_1,j)$ and $\sigma(i_2,j)$ for $j>1$ are
disjoint, i.e., the images of any element are distinct under these permutations;

(iii) similar condition holds for $\sigma(i,j_1)$ and $\sigma(i,j_2)$ for $i>1$.

Both yet unproved statements of the theorem follow.
                            \qed

\medskip
It is intersting to notice that in terms of line arrangements the above
full-graph connected
components of $\Gamma$ correspond to sets of lines such that any point from $\X$
lies on exactly one line of the set. We will discuss the existence of such
arrangements in 6.3.

We obtain another series of representable matrices $Q$
using the affine matrices corresponded to
bushes.
Let $\Gamma$ be the graph of three connected components, each one being 
the bush on $r+1$ vertices. Let the restriction of the labeling vector $m$
to each connected component be $(1,r,\ldots,r)$ where 1 labels the root of the
bush. Then $Q(m,\Gamma)$ has three affine components and is realizable. As a
matroid one can take that of
the complex
arrangement of lines of type $\A_{r,1,3}$ given by
$xyz(x^r-y^r)(x^r-z^r)(y^r-z^r)$. It follows immediately from \cite{CS}, Remark
6.3.

\medskip
In the rest of this subsection we classify all the realizable matrices
with several affine components such that at least three of the components
are Cartan matrices.
Our argument is
straightforward: given an indecomposable matrix we consider all possible
realizations of it and list all transversal columns for each realization.
Going through this list we can see what other components can be realized jointly
with the initial matrix. In order to reduce amount of case by case inspections the
following simple observation is useful.

\begin{lemma}
\label{subgraph}
Let $\Gamma$ be a graph on vertices $\{v_1,\ldots,v_n\}$,
 $m=(m_1,\ldots,m_n)$ a labeling of its vertices and $Q=Q(m,\Gamma)$.
Fix $K\in\Pi$ and $i\in K$.
Suppose $J$ is a realization of $Q$ and $J'$ is obtained from $J$ by deleting
the $i$th column and all the rows whose restrictions to $K\setminus\{i\}$ 
are 0. Then $J'$ is a realization of $Q'=Q(m',\Gamma')$ where
$\Gamma'$ is the subgraph of $\Gamma$ on
$\{v_1,\ldots,\hat v_i,\ldots,v_n\}$ and $m'_j\leq m_j$ for all $j\not=i$ 
with equality on $K\setminus\{i\}$.
\end{lemma}

Since we consider Cartan matrices the labels $m$ have $m_i=2$ for all $i$.
Notice that in terms of arrangements of lines, this means that we consider only
systems of multiple
points such that any line passing through them contains precisely 3 of them.
 Under this condition, the indecomposable affine matrices are
Cartan matrices of extended Dynkin diagrams with singular edges, i.e., the types
$A_{\ell}^{(1)}$ ($\ell\geq 2$), $D_{\ell}^{(1)}$ ($\ell\geq 4$) and
$E_{\ell}^{(1)}$ ($\ell=6,7,8$).
Notice that all these graphs except $A_2^{(1)}$ contain $\Gamma_0=A_3$ as a subgrpah.
Thus the first graph to analyze is $\Gamma_0$. 

The following facts about $\Gamma_0$ can be obtained by simple inspection.

(i) The Cartan matrix of $\Gamma_0$ has the unique (up to permutations)
realization by
$$J_0=\left (\matrix{1&0&1\cr
                   1&0&0\cr
                   1&0&0\cr
                   0&1&0\cr
                   0&1&0\cr
                   0&1&0\cr
                   0&0&1\cr
                   0&0&1\cr}\right).$$
(ii) There are precisely two (up to permutations of columns and rows) maximal sets
of columns transversal to the collection of columns of $J_0$ such that every two of
them intersect at no more than one element. Written as matrices
they are
$$T_1=\left(\matrix{1&1&1&0&0&0&0\cr
                    0&0&0&1&1&0&0\cr
                    0&0&0&0&0&1&1\cr
                    1&0&0&1&0&0&1\cr
                    0&1&0&0&1&1&0\cr
                    0&0&1&0&0&0&0\cr
                    0&0&0&1&0&1&0\cr
                    0&0&0&0&1&0&1\cr}\right),
  T_2=\left(\matrix{1&1&1&0&0&0&0\cr
                    0&0&0&1&1&0&0\cr
                    0&0&0&0&0&1&1\cr
                    1&0&0&1&0&0&0\cr
                    0&1&0&0&1&1&0\cr
                    0&0&1&0&0&0&1\cr
                    0&0&0&1&0&1&0\cr
                    0&0&0&0&1&0&1\cr}\right).$$
The respective symmetric matrices $R_i=T_i^tT_i-E$ ($i=1,2$) are
$$R_1=\left(\matrix{1&0&0&0&-1&-1&0\cr
                    0&1&0&-1&0&0&-1\cr
                    0&0&1&-1&-1&-1&-1\cr
                    0&-1&-1&2&0&0&0\cr
                    -1&0&-1&0&2&0&0\cr
                    -1&0&-1&0&0&2&0\cr
                    0&-1&-1&0&0&0&2\cr}\right),$$
  $$R_2=\left(\matrix{1&0&0&0&-1&-1&-1\cr
                    0&1&0&-1&0&0&-1\cr
                    0&0&1&-1&-1&-1&0\cr
                    0&-1&-1&2&0&0&-1\cr
                    -1&0&-1&0&2&0&0\cr
                    -1&0&-1&0&0&2&0\cr
                    -1&-1&0&-1&0&0&2\cr}\right).$$

These matrices correspond to the labeled graphs $\Gamma_1$ and $\Gamma_2$.

\PSbox{fig1.pstex}{10cm}{10cm}

$$\matrix{&&&&&& \Gamma_1 &&&&&&&&&&&&&&&&&&&&&&
 \Gamma_2&&&&& }$$

\bigskip
(iii) Now suppose that a graph $\Gamma$ has (at least three)
extended Dynkin diagrams with single edges as its
connected components and at least one of them, say $\Gamma'$,
 contains $A_3$ as a subgraph (i.e.,
this component is not $A_2^{(1)}$). 
Suppose the respective Cartan matrix $Q=Q(m,\Gamma)$ ($m=(2,\ldots,2)$)
is realizable by $J$. Denote by $J'$ the submatrix of $J$ that consists of 
the columns of $J$ corresponding to vertices of $\Gamma'$ and by $J''$ the matrix
of the other columns. 
We can apply Lemma \ref{subgraph}.
It implies that $J''$ can be obtained from either of $J_1$ and $J_2$
by adding three new rows that form the identity block with the first
three rows of $J_i$ and have only 0 in the other positions and then
omitting some columns (maybe none). The addition of three new rows
does not change $\Gamma_i$, just makes the labeling
$(2,2,\ldots,2)$. Thus one should be able to delete some vertices from
one of these graphs and obtain a graph having at least two extended
Dynkin diagrams as its connected components. It is quite obviously
impossible and we get a contradiction. Thus $Q$ cannot be
realized. 

(iv) It is left to consider $Q$ whose each component is the Cartan matrix
$A_2^{(1)}$. This case is covered by Theorem \ref{full graphs}. Moreover since 
the permutation group on $\{1,2,3\}$ 
has only two elements without a fixed point there is the unique realization of $Q$
(up to permutations of rows and columns). The maximum number of components $Q$ can
have is 4 and if it has 4 affine components it is realized by the matrix $J_C$
$$J_C=\left(\matrix{1&0&0&1&0&0&1&0&0&1&0&0\cr
                  1&0&0&0&1&0&0&1&0&0&1&0\cr
                  1&0&0&0&0&1&0&0&1&0&0&1\cr
                  0&1&0&1&0&0&0&1&0&0&0&1\cr
                  0&1&0&0&1&0&0&0&1&1&0&0\cr
                  0&1&0&0&0&1&1&0&0&0&1&0\cr
                  0&0&1&1&0&0&0&0&1&0&1&0\cr
                  0&0&1&0&1&0&1&0&0&0&0&1\cr
                  0&0&1&0&0&1&0&1&0&1&0&0\cr}\right).$$
If it has three affine components then it is realized by either 
9 or 10 or 11 first columns of $J_C$.

Summing up, we obtain the following result. 

\begin{proposition}
\label{cartan}
Let $Q$ have at least three affine components that are Cartan matrices. Then
$Q$ is realizable if and only if $Q$ has at most four components such that
each
affine component is $2I_3$ (where $I_3$ is the 3$\times$3 identity matrix)
and if a finite component exists then it is a principle submatix of $2I_3$.
\end{proposition}

\medskip
Similarly it is possible to study the case where $m_i\leq 2$ for every $i$. The
only extra possibility here is given by three components each being the
affine matrix corresponding to the bush with 2 leaves.

\medskip
\noindent{\bf 6.3. Representations of matroids.}

\medskip
\noindent If a matrix $Q$ is realizable then we can ask if the respective matroid
can be represented by an arrangement $\A$ of lines in a projective plane
over a given field. If it is possible then the set $\X$ of 
all multiple points of $\A$ define 
component $V(\X)$ of $R_1$ with the associated matrix $Q$. 

Here we give the positive answers about the representability of certain
matroids from
this section. 

For $Q$ from Theorem
\ref{full graphs} 
representations over $\cn$ exist for $n=2$ and $n=3$. For $n=2$ it
is the braid arrangement given by $xyz(x-y)(x-z)(y-z)$ and for $n=3$ the
Hessian arrangement of 12 lines containing 9 flexes of a non-singular cubic
(e.g., given by $xyz\prod_{i,j=0,1,2}(x+\zeta^iy+\zeta^jz)$ where
$\zeta={\rm exp}(2\pi\sqrt{-1}/3)$). 
In Section 7, we will prove that for $n>3$ these matroids cannot be represented
over $\cn$.

 The arrangement of 9 lines projectively dual to the Hessian arrangement (that is
a representation of the plane of order 3)
gives a different example. All the 12 points
of this arrangement do not define any component of $R_1$ since all components of the
respective matrix $Q$ are finite. Any subset of 9 points such that the left out 3
points have all 9 lines passing through them defines a component of $R_1$ 
of the above type with $n=\ell=3$.

Let us remark about matrices
$Q$ from Theorem \ref{full graphs} that they completely
cover the case of arrangements of lines with only triple multiple points and
$R_1\not=0$. Indeed if any row of matrix $J$ has exactly three 1's then $Q$ has
at most three components. Since $R_1\not=0$ the number of components is exactly
three and each one of them corresponds to a full graph.
It follows from the above theorem that there are $3n$ lines (for some $n$)
partitioned in 3 groups of equal size $n$ such that lines in each group are in
general position and each multiple point is on a line from each group.
In section 7, we will show that these arrangements (equivalently their duals)
do not exist for $n>5$.

Another series of representable (by definition) matroids come from realizations
of matrices having components corresponding to bushes. The respective
arrangements are
of types $\A_{r,1,3}$.

Finally Proposition \ref{cartan} gives the following result for arrangements.

\begin{theorem}
\label{3-lines}
Let $\A$ be a projective arrangement of lines with a set $\X$ of points from $L$
such that $V(\X)$ is a component of $R_1$ of positive dimension. Let any line from
$L$ having a point from $\X$ have exactly three of them. Then 
$L(2)$ for this arrangement is given by the incidence
matrix that is a submatrix of $J_C$ from subsection 6.2 that 
includes the first 9 of its columns.
\end{theorem}

\bigskip 
\section{Direct sum decompositions and pencils.}

\bigskip 

In the case where $\A$ is an arrangement over $\cn$ the 
matrix $Q=Q(\X)$ assigned to a subset $\X$ of $L'(2)$
has a natural geometric  interpretation as 
follows. Let us consider the arrangement of lines induced by $\A$
in a generic plane
$P$ in ${{\bf P}(\cn)}^{\ell-1}$. Let $P'$ be the blow up of 
$P$ at the intersection points of the images of $X\in \X$ with $P$.
Let $\Sigma$ be a free abelian  group generators of which correspond to the 
proper preimages in $P'$ 
of  intersections  $P$ with the   hyperplanes from
$\cal A$. 

\begin{proposition}
\label{intersection matrix} 
Intersection form on $P'$ induces a symmetric bilinear form on $\Sigma$
which has $-Q$ as its matrix.
\end{proposition}

\proof 
The entries of the bilinear form induced on $\Sigma$ by  
the intersection form of $P'$ are the following. 
If the lines in $P$ intersect at a point from $\X$
then the entry corresponding to these lines is $0$, if the lines
intersect at a point outside of $\X$ then the entry is $1$ and the
diagonal entry is $1-m$ where $m$ is the number of points from $\X$ 
on this line (all this since a blow up at a point decreases the
intersection index by 1). Comparing with the definition of $Q$
in section 2 we conclude the proof. \qed

\medskip 

Now first part of 
Proposition 2.2 (ii) from this point of view is a consequence of the 
Hodge index theorem claiming that the signature  of the intersection 
form on $P'$ is $(1,rk H_2(P')-1)$
 
\medskip
Next let us assume that $\dim R_1 >0$ and the corresponding
component $S$ of the characteristic variety (i.e. the image of a
component of $R_1$ under the exponential map) is essential 
(cf. \cite{L}). The latter means that $\cal A$ does not  
have a subarrangement ${\cal A}^{\prime}$ such that for some component $S'$ of 
the characteristic variety of ${\cal A}^{\prime}$ one has $S=i^*(S')$ where
the map $i^*$ of the character groups  
$i^*: H^1({\bf P}^2(\cn)\setminus{\cal A}^{\prime},{\bf C}^*) \rightarrow 
 H^1({\bf P}^2(\cn)\setminus{\cal A},{\bf C}^*)$ is induced by inclusion.
Recall also that $S$ corresponds to a choice of a subset $\X 
\subset L'(2)$ and to emphasize this we sometimes write
$S(\X)$ for $S$.

Let $\pi: P\setminus{\cal A} \rightarrow {\bf P}^1\setminus\{p_1,\ldots,p_m\}$ 
be the map corresponding to the component $S$ (cf. section 5).
The map $\pi$ defines a rational map $\tilde \pi:  P \rightarrow {\bf P}^1$ 
regular outside of a set of codimension 2 (cf. \cite{Mumford}), 
i.e. a finite set of points. Since $\pi$ is regular outside of 
${\cal A}$ we have
 $$\bigcup \tilde \pi^{-1} (p_i) \subset {\cal A}. \eqno (7.1)$$
The pencil of curves consisting of the closures of fibres of $\tilde
\pi$ does not have common components. Indeed, if we assume that $\ell$ is such an
irreducible  component, then from (7.1) $\ell$ must  
be a line and $S$ is the image of a component in the arrangement 
obtained from $\cal A$ by deleting $\ell$ in contradiction with
assumption that $S$ is essential. This implies that the pull back 
of the pencil corresponding to $\tilde \pi$ has no base points on 
$P'$ (i.e. the blow up of $P$ at $\X$). Let 
$\pi': P' \rightarrow {\bf P}^1$ be the corresponding map. 
It follows from Bertini's theorem  that a generic fibre of $\pi'$ is non
singular. Now let us consider a free abelian group $D_i$ generated by 
the classes of lines from $\cal A$ which belong 
to the fibre $\pi'^{-1}(p_i)$. Clearly subgroups $D_i$ and $D_j$ 
of $\Sigma$ are orthogonal
with respect to the intersection form. On the other hand, since neither
of the exceptional curves in $P'$ belongs to a fibre of $\pi'$ and by 
Zariski's connectedness theorem the fibres of $\pi'$ are connected
(cf. \cite{Mumford} (3.24)), it follows that the intersection form 
is irreducible on $D_i$. 

Next recall (\cite{monodromy}, exp. X) that for each fibre of a fibration 
of a surface there is unique linear combination  $z_0$ of components 
of this fibre with positive integer coefficients such that 

\par a) for any combination $z$ of the components of this fibre 
$(z,z) \le 0$ with equality if and only if $z$ is an integer multiple of 
$z_0$. 

\par b) $(z,z_0)=0$ for any $z$ supported on the fibre. 

A consequence of this is that the intersection form on each $D_i$ is
affine. 

We can summarize this as  follows.

\begin{proposition}
\label{fibres affine} Let $\cal A$ be an arrangement such that 
$\dim R_1 >0$ and $S=S(\X)$ an essential component of a characteristic 
variety containing 1.
Then there is a one-to-one correspondence between the blocks of matrix 
$Q(\X)$ from section 2 and the fibres of the pencil corresponding to the  
component $S(\X)$ .
\end{proposition}

\medskip Existence of a pencil for which given arrangement is a union 
of several elements imposes severe restrictions on the arrangement. 
For example one can show the following (cf. 6.3).

\begin{proposition}
\label{classification} Let $\cal A$ be an arrangement of $n \cdot k$ ($n
\ge 2$) lines for which: 
\par a) the corresponding matrix $Q$ is sum of $k$ equal blocks so that 
the lines corresponding to each block are in general position; 
\par b) any multiple point contains a line from each group;
\par c) each line contains $n$ multiple points of the arrangement. 
\par Then $k \le 5$.
\end{proposition} 

\proof It follows from c) that each group is affine and hence $\dim
R_1>0$. Hence there is a pencil such that the arrangement forms a union 
of its singular fibres.
Let us calculate the euler characteristic of $P'$, i.e. the plane blown
up at the multiple points of the arrangement. According to c) and a)
there are $n^2$ multiple points and hence the euler characteristic
of $P'$ is $E_1=3+n^2$. On the other hand, since $P'$ is fibered over
${\bf P}^1$ with $k$ fibers consisting of $n$ lines in general
position (which has the euler characteristic equal to $2n-{{n(n+1)} \over 2}$)
and with generic fibre a non singular plane curve of degree $n$ (which 
has the euler characteristic $3n-n^2$), we see that  $e(P') \ge E_2$
where $E_2=(2-k)) \cdot (3n-n^2)+k \cdot (2n- {{n(n+1)} \over
2})$ (inequality is due to possible presence of other singular fibres
and the semicontinuity of euler characteristic). But $E_2 >E_1$ for $k \ge
6$ which yields (i). Also in the case $k=n+1$ we see that 
$E_1=E_2$ if and only if  $n=2,3$, i.e. in these cases we may and indeed 
do have a fibration with 
no degenerate fibres other than those composed of the lines of
arrangement.
\qed

\medskip 

Similar argument allows to obtain a restriction on arrangements having
a given number of multiple point on each line which bounds the number of
arrangements having large $\dim R_1$. 

\begin{proposition} 
\label{asymptotic} Let ${\cal A}(r,k)$ be a class of arrangements with 
each line having at most $k$ points and such that $dim R_1 =r$. 
Then there is a function $F(r,k) \le 2k+1$ 
such that the number of lines in each of $r+1$ groups 
in which the arrangement splits 
does not exceed $F(r,k)$. In particular the number of the lines in the 
arrangement is at most $(r+1)F(r,k) < (r+1)(2k+1)$. 
\end{proposition}

\proof The argument is similar to the one used in the previous 
proposition. Let $d$ be the maximum of the numbers of lines in each 
group of lines forming a singular fibre of the pencil. Then the euler
characteristic of the total space of the fibration is 
greater than $E_2=(2-r)(3d-d^2)+r(2d-d(d-1)/2)$ (use semicontinuity of 
euler characteristic). On the other hand the euler characteristic 
of the blow up is greater than or equal to $E_1=3+d \cdot k \cdot r$.
One has $E_2 > E_1$ and hence a contradiction unless $E_2-E_1$ has 
a positive root $F(r,k)$ which is the bound on $d$. Direct calculation
shows that when $r \rightarrow \infty$ the positive root of
$E_2-E_1$ is an increasing function of $r$ and its limit is $2k+1$. \qed

\end{document}